\documentclass[10pt]{article}
\usepackage{times,latexsym,amsfonts}
\usepackage[dvips]{graphics}

\def\bmatrix#1{\left[\matrix{#1}\right]}

\newcommand{\QED}{\hspace*{\fill}\mbox{\footnotesize$\Box$}\par\bigskip}

\def\kmtheoremdef#1#2#3#4{
  \newtheorem{#1x}{#2}[section]
  \newenvironment{#1}{\begin{#1x}#3\relax}{#4\end{#1x}}}

\kmtheoremdef{assumption}{Assumption}{\rm}{}
\kmtheoremdef{corollary}{Corollary}{\it}{}
\kmtheoremdef{definition}{Definition}{\rm}{}
\kmtheoremdef{example}{Example}{\rm}{}
\kmtheoremdef{lemma}{Lemma}{\it}{}
\kmtheoremdef{note}{Note}{\rm}{}
\kmtheoremdef{proposition}{Proposition}{\it}{}
\kmtheoremdef{theorem}{Theorem}{\it}{}

\newenvironment{proof}{\par\noindent{\textbf{Proof}}.\ }{\QED}

\newcommand{\figureherecaption}[1]{Figure\,\thefigure: #1}
\newenvironment{figurehere}{
  \begin{center}%
  \refstepcounter{figure}%
}{%
  \end{center}%
}

\renewcommand{\Bigl}[1]{\mbox{\large\boldmath$\bigl#1$}}
\renewcommand{\Bigr}[1]{\mbox{\large\boldmath$\bigr#1$}}

\begin{document}
\begin{center}
\textbf{\Large
  A~Parameterization of Stabilizing Controllers
             over Commutative Rings
}

\bigskip
\bigskip
\begin{large}
  Kazuyoshi MORI$\dag\ddag$
\end{large}

$\dag$ Department of Electrical Engineering, Faculty of Engineering,\\
       Tohoku University, Sendai 980-8579, JAPAN\\
       (\texttt{Kazuyoshi.MORI@IEEE.ORG})
\end{center}

\renewcommand{\thefootnote}{\fnsymbol{footnote}} \footnotetext[3]{
Part of this work was done at 
Institut de Recherche en Cybern\'etique de Nantes, France
and the author thanks this institution for its hospitality and
support.}
\renewcommand{\thefootnote}{\arabic{footnote}}

\section{Introduction}\label{S:Introduction}
In this paper we are concerned with the factorization approach to
control systems, which has the advantage that it embraces, within
a~single framework, numerous linear systems such as continuous-time as
well as discrete-time systems, lumped as well as distributed systems,
$1$-D as well as~$n$-D systems, etc.\cite{bib:vidyasagar82a}.  The
factorization approach was patterned after Desoer \emph{et
al.}\cite{bib:desoer80a} and Vidyasagar \emph{et
al.}\cite{bib:vidyasagar82a}.  In this approach, when problems such as
feedback stabilization are studied, one can {focus} on the key aspects
of the problem under study rather than be distracted by the special
features of a~particular class of linear systems.  A~transfer matrix
of this approach is considered as the ratio of {two} stable causal
transfer matrices.  For a~long time, the theory of the factorization
approach had been founded on the coprime factorizability of transfer
matrices, which is satisfied by transfer matrices over the principal
ideal domains or the Bezout domains.

Sule in~\cite{bib:sule94a,bib:sule98a} has presented a~theory of the
feedback stabilization of strictly causal plants for multi-input
multi-output transfer matrices over commutative rings with some
restrictions.  This approach to the stabilization theory is called
``coordinate-free approach'' in the sense that the coprime
factorizability of transfer matrices is not required.  Recently, Mori
and Abe in~\cite{bib:mori97a,bib:mori98bs2} have generalized his theory
over commutative rings under the assumption that the plant is causal.
They have introduced the notion of the generalized elementary factor,
which is a~generalization of the elementary factor introduced by
Sule\cite{bib:sule94a}, and have given the necessary and sufficient
condition of the feedback stabilizability.

Since the stabilizing controllers are not unique in general, the
choice of the stabilizing controllers is important for the resulting
closed loop.  In the classical case, that is, in the case where there
exist the right-/left-coprime factorizations of the given plant, the
stabilizing controllers can be parameterized by the method called
``Youla-Ku\v{c}era
parameterization''\cite{bib:vidyasagar82a,bib:desoer80a,bib:raman84a,bib:youla76a}.
However, it is known that there exist models in which some
stabilizable transfer matrices do not have their right-/left-coprime
factorizations\cite{bib:anantharam85a}.  In such models, we cannot
employ the Youla-Ku\v{c}era parameterization directly.  In this paper
we will give a~parameterization of the stabilizing controllers over
commutative rings by using the results given by
Sule\cite{bib:sule94a}, and Mori and
Abe\cite{bib:mori97a,bib:mori98bs2}.

Here we briefly outline how the parameterization, which is different
from the Youla-Ku\v{c}era parameterization, will be obtained.  Let
$H(P,C)$ denote the transfer matrix of the standard feedback system
defined as
\[
  H(P,C) =
   \bmatrix{
    (E+PC)^{-1}  &  -P(E+CP)^{-1} \cr
    C(E+PC)^{-1} & (E+CP)^{-1}},
\]
where~$P$ and~$C$ are plant and controller, and~$E$ the identity
matrix.  We consider the set~${\cal H}$ of $H(P,C)$'s with all stabilizing
controllers~$C$ rather than the set of all stabilizing controllers
itself.  We will characterize this~${\cal H}$ by one parameter matrix.  Then
using it, we will obtain the parameterization of the stabilizing
controllers.

The paper is organized as follows.  After this introduction, we begin
on the preliminary in~\S\,\ref{S:Preliminary}, in which we give
mathematical preliminaries, set up the feedback stabilization problem
and present the previous results.  To obtain the set~${\cal H}$ above we
will use both of right-/left-coprime factorizations over the ring of
fractions of the set of the stable causal transfer functions.  In
order to establish the existence of such right-/left-coprime
factorizations, we will present, in~\S\,\ref{S:Stab}, the one-to-one
correspondence between the sets of the radicals of the generalized
elementary factors of the plant and its transposed plant.
Section\,\ref{S:01.Feb.99.173244} is the main part of this paper, in
which a~parameterization of the stabilizing controllers is presented.
In~\S\,\ref{S:Example} we will consider the multidimensional system
with structural stability as an example and present the
parameterization of its stabilizing controllers.  Our method will give
a solution of an open problem in\,\cite{bib:lin98a} about the
parameterization of the stabilizing controllers for the
multidimensional system with structural stability.

\section{Preliminaries}
\label{S:Preliminary}
In the following we begin by introducing the notations of commutative
rings, matrices, and modules used in this paper.  Then we give the
formulation of the feedback stabilization problem and the previous
results.  We also review the construction method of a~stabilizing
controller presented in~\cite{bib:mori98bs2}.

\subsection{Notations}\label{SS:Notations}
\paragraph{Commutative Rings}
In this paper we consider that any commutative ring has the
identity~$1$ different from zero.  Let~${\cal R}$ denote a~(unspecified)
commutative ring.  The total ring of fractions of~${\cal R}$ is denoted by~$
{\cal F}({\cal R})$.

We will consider that \emph{the set of stable causal transfer
functions} is a~commutative ring denoted by~${\cal A}$.  Further we will use
the following three kinds of ring of fractions.
The first one appears as the total ring of fractions of~${\cal A}$, which is
denoted by~${\cal F}({\cal A})$ or simply by~${\cal F}$; that is, ${\cal F}=\{ n/d\,|\, n,d\in{\cal A},~
\mbox{$d$ is a~nonzerodivisor}\}$.  This will be considered as
\emph{the set of all possible transfer functions}.
The second one is associated with the set of powers of a~nonzero
element of~${\cal A}$.  Let~$f$ denote a~nonzero element of~${\cal A}$.  Given
a~set $S_f=\{ 1,f,f^2,\ldots\}$, which is a~multiplicative subset of~${\cal A}$,
we denote by~${\cal A}_f$ the ring of fractions of~${\cal A}$ with respect to the
multiplicative subset~$S_f$; that is, ${\cal A}_f=\{ n/d\,|\, n\in{\cal A},~d\in S_f\}$.
It should be noted that, in the case where~$f$ is a~zerodivisor, even
if the equality $a/1=b/1$ with $a,b\in{\cal A}$ holds over~${\cal A}_f$, we cannot
say in general that $a=b$ over~${\cal A}$; alternatively, $a=b+z$ over~${\cal A}$
holds with some zerodivisor~$z$ of~${\cal A}$ such that $zf=0$.
The last one is the total ring of fractions of~${\cal A}_f$, which is
denoted by ${\cal F}({\cal A}_f)$; that is, ${\cal F}({\cal A}_f)=\{ n/d\,|\, n,d\in{\cal A}_f,~\mbox{$d$
is a~nonzerodivisor of ${\cal A}_f$}\}$.  If~$f$ is a~nonzerodivisor of~${\cal A}
$, ${\cal F}({\cal A}_f) $ coincides with the total ring of fractions of~${\cal A}$.
Otherwise, they do not coincide.  The reader is referred to~Chapter\,3
of~\cite{bib:atiyah69a} for the ring of fractions.

\paragraph{Matrices}
The set of matrices over~${\cal R}$ of size $x\times y$ is denoted by ${\cal R}^{x\times
y}$.  Further, the set of square matrices over~${\cal R}$ of size~$x$ is
denoted by~$ ({\cal R})_x$.  The identity and the zero matrices are denoted
by~$E_x$ and $O_{x\times y}$, respectively, if the sizes are required,
otherwise they are denoted by~$E$ and~$O$.

Matrix~$A$ over~${\cal R}$ is said to be \emph{nonsingular}
$\Bigl($\emph{singular}$\Bigr)$ \emph{over~${\cal R}$} if the determinant of
the matrix~$A$ is a~nonzerodivisor $\Bigl($a zerodivisor$\Bigr)$ of~$
{\cal R}$.  Matrices~$A$ and~$B$ over~${\cal R}$ are \emph{right-}
\emph{$\Bigl($left-$\Bigr)$coprime over~${\cal R}$} if there exist
matrices~$X$ and~$Y$ over~${\cal R}$ such that
  $XA+YB=E$ $\Bigl(AX+BY=E\Bigr)$ 
holds.  Note that, in the sense of the above definition, two matrices
which have no common right-$\Bigl($left-$\Bigr)$divisors except
invertible matrices may not be right-$\Bigl($left-$\Bigr)$coprime
over~${\cal R}$.  (For example, two matrices $\bmatrix{z_1}$ and
$\bmatrix{z_2}$ of size~$1\times 1$ over the bivariate polynomial ring
$\mathbb{R}[z_1,z_2]$ over the real field $\mathbb{R}$ are neither
right- nor left-coprime over $\mathbb{R}[z_1,z_2]$ in our setting.)
Further, an ordered~pair $(N,D)$ of matrices~$N$ and~$D$ is said to be
a \emph{right-coprime factorization over~${\cal R}$} of~$P$ if (i)~$D$ is
nonsingular over~${\cal R}$, (ii) $P=ND^{-1}$ over ${\cal F}({\cal R})$, and (iii)~$N$
and~$D$ are right-coprime over~${\cal R}$.  As the parallel notion, the
\emph{left-coprime factorization over~${\cal R}$} of~$P$ is defined
analogously.  That is, an ordered~pair $(\widetilde{D},\widetilde{N})$
of matrices~$\widetilde{N}$ and~$\widetilde{D}$ is said to be a
\emph{left-coprime factorization over~${\cal R}$} of~$P$ if (i)
$\widetilde{D}$ is nonsingular over~${\cal R}$, (ii)
$P=\widetilde{D}^{-1}\widetilde{N}$ over ${\cal F}({\cal R})$, and
(iii)~$\widetilde{N}$ and~$\widetilde{D}$ are left-coprime over~${\cal R}$.
Note that the order of the ``denominator'' and ``numerator'' matrices
is interchanged in the latter case.  This is to reinforce the point
that if $(N,D)$ is a~right-coprime factorization over~${\cal R}$ of~$P$,
then $P=ND^{-1}$, whereas if $(\widetilde{D},\widetilde{N})$ is
a~left-coprime factorization over~${\cal R}$ of~$P$, then
$P=\widetilde{D}^{-1}\widetilde{N}$ according
to~\cite{bib:vidyasagar85a}.  For short, we may omit ``over~${\cal R}$''
when ${\cal R}={\cal A}$, and ``right'' and ``left'' when the size of matrix is $1
\times 1$.

\paragraph{Modules}
Let $M_r(X)$ \,$\Bigl(M_c(X)\Bigr)$ denote the~${\cal R}$-module generated
by rows \,$\Bigl($columns$\Bigr)$ of a~matrix~$X$ over~${\cal R}$.  Let
$X=AB^{-1}=\widetilde{B}^{-1}\widetilde{A}$ be a~matrix over~${\cal F}({\cal R})$,
where~$A$,~$B$,~$\widetilde{A}$,~$\widetilde{B}$ are matrices over~${\cal R}
$.  It is known that $M_r(\bmatrix{A^t & B^t}^t)$
\,\,$\Bigl(M_c(\bmatrix{\widetilde{A} & \widetilde{B}})\Bigr)$ is
unique up to an isomorphism with respect to any choice of fractions
$AB^{-1}$ of~$X$ \,$\Bigl(\widetilde{B}^{-1}\widetilde{A}$ of
$X\Bigr)$ (Lemma\,2.1 of\,\cite{bib:mori98bs2}).  Therefore, for
a~matrix~$X$ over~${\cal R}$, we denote by~${\cal T}_{X,{\cal R}}$ and~${\cal W}_{X,{\cal R}}$ the
modules $M_r(\bmatrix{A^t & B^t}^t)$ and $M_c(\bmatrix{\widetilde{A} &
\widetilde{B}})$, respectively.  

\subsection{Feedback Stabilization Problem}\label{SS:FSProblem}
The stabilization problem considered in this paper follows that of
Sule in~\cite{bib:sule94a}, and Mori and Abe
in~\cite{bib:mori97a,bib:mori98bs2}, who consider the feedback system~$
\Sigma$~\cite[Ch.5, Figure\,5.1]{bib:vidyasagar85a} as in
Figure\,\ref{Fig:FeedbackSystem}.
\begin{figurehere}
\begin{center}
  \resizebox{0.4\textwidth}{!}{\includegraphics{feedbacksystem.eps}}\\
  \figureherecaption{Feedback system~$\Sigma$.}\label{Fig:FeedbackSystem}
\end{center}
\end{figurehere}
For further details the reader is referred
to~\cite{bib:vidyasagar85a}.  Throughout the paper, the plant we
consider has~$m$ inputs and~$n$ outputs, and its transfer matrix,
which is also called a \emph{plant} itself simply, is denoted by~$P$
and belongs to~${\cal F}^{n\times m}$.  We can always represent~$P$ in the form
of a~fraction $P=ND^{-1}$
$\Bigl(P=\widetilde{D}^{-1}\widetilde{N}\Bigr)$, where $N\in{\cal A}^{n\times m}$
\,$\Bigl(\widetilde{N}\in{\cal A}^{n\times m}\Bigr)$ and $D\in({\cal A})_m$
\,$\Bigl(\widetilde{D}\in({\cal A})_n\Bigr)$ with nonsingular~$D$
\,$\Bigl(\widetilde{D}\Bigr)$.  
\begin{definition}
For $P\in{\cal F}^{n\times m}$ and $C\in{\cal F}^{m\times n}$, a~matrix $H(P,C)\in({\cal F})_{m+n}$ is
defined as
\begin{equation}\label{E:H(P,C)}
  H(P,C) =
   \bmatrix{
    (E_n+PC)^{-1}  &  -P(E_m+CP)^{-1} \cr
    C(E_n+PC)^{-1} & (E_m+CP)^{-1}}
\end{equation}
provided that $\det(E_n+PC)$ is a~nonzerodivisor of~${\cal A}$.  This
$H(P,C)$ is the transfer matrix from $\bmatrix{u_1^t & u_2^t}^t$ to
$\bmatrix{e_1^t & e_2^t}^t$ of the feedback system~$\Sigma$.  If (i)
$\det(E_n+PC)$ is a~nonzerodivisor of~${\cal A}$ and (ii) $H(P,C)\in({\cal A})_
{m+n}$, then we say that the plant~$P$ is \emph{stabilizable}, $P$ is
\emph{stabilized} by~$C$, and~$C$ is a \emph{stabilizing controller}
of~$P$.
\end{definition}

Since the transfer matrix $H(P,C)$ of the stable causal feedback
system has all entries in~${\cal A}$, we call the above notion~\emph{${\cal A}
$-stabilizability}.  One can further introduce the notion of~\emph{${\cal A}_
f$-stabilizability} as follows.
\begin{definition}
Let~$f$ be a~nonzero element of~${\cal A}$.  If (i) $\det(E_n+PC)$ is a~nonzerodivisor of~${\cal A}_f$ and (ii) $H(P,C)\in({\cal A}_ f)_{m+n}$, then we say
that the plant~$P$ is \emph{${\cal A}_ f$-stabilizable}, $P$ is \emph{${\cal A}_
f$-stabilized} by~$C$, and~$C$ is an \emph{${\cal A}_f$-stabilizing
controller} of~$P$.
\end{definition}

The causality of transfer functions is an important physical
constraint.
We employ, in this paper, the definition of the causality from
Vidyasagar~\emph{et al.}\cite[Definition\,3.1]{bib:vidyasagar82a}.
\begin{definition}\label{D:Causal}
Let~${\cal Z}$ be a~prime ideal of~${\cal A}$, with ${\cal Z}\neq{\cal A}$, including all
zerodivisors.  Define the subsets~${\cal P}$ and ${\cal P}_{\textrm{s}}$ of~${\cal F}$
as follows:
\begin{eqnarray*}
  {\cal P}&=&\{ n/d\in{\cal F}\,|\, n\in{\cal A},~ d\in{\cal A}\backslash{\cal Z}\},\\
  {\cal P}_{\textrm{s}}&=&\{ n/d\in{\cal F}\,|\, n\in{\cal Z},~ d\in{\cal A}\backslash{\cal Z}\}.
\end{eqnarray*}
Then every transfer function in~${\cal P}$ $\Bigl({\cal P}_{\textrm{s}}\Bigr)$ is
called \emph{causal} $\Bigl($\emph{strictly causal}$\Bigr)$.
Analogously, if every entry of a~transfer matrix~$F$ is in~${\cal P}$
$\Bigl({\cal P}_ {\textrm{s}}\Bigr)$, the transfer matrix~$F$ is called
\emph{causal} $\Bigl($\emph{strictly causal}$\Bigr)$.  A~matrix over~$
{\cal A}$ is said to be \emph{${\cal Z}$-nonsingular} if the determinant is in ${\cal A}\backslash
{\cal Z}$, and \emph{${\cal Z}$-singular} otherwise.
\end{definition}

\subsection{Previous Results}
\label{S:PR}
To state precisely the previous results of the stabilizability, we
recall the notion of the generalized elementary factors, which is
a~generalization of the elementary factor given by
Sule\cite{bib:sule94a}.  Originally, the elementary factor has been
defined over unique factorization domains.  Mori and Abe have enlarged
this concept for commutative rings\cite{bib:mori97a,bib:mori98bs2}.

\begin{definition}\label{D:GEF}
\textup{(Generalized Elementary Factors, Definition\,3.1
of~\cite{bib:mori98bs2})~} Let ${\cal I}
$
be the set of all sets of~$m$ distinct integers between~$1$ and $m+n$.
Let~$I$ be an element of~${\cal I}$ and $i_1,\ldots, i_m$ be elements of~$I$
with ascending order, that is, $i_a<i_b$ if $a<b$.  We will
use~elements of~${\cal I}$ as suffices as well as sets.  Denote by~$T$ the
matrix $\bmatrix{N^t&D^t}^t$ over~${\cal A}$ with $P=ND^{-1}$.  Let $\Delta_I\in
{\cal A}^{m\times(m+n)}$ denote the matrix whose $(k,i_k)$-entry is~$1$ for $i_k
\in I$ and zero otherwise.  For each $I\in{\cal I}$, an ideal $\Lambda_ {P\!I}$
over~${\cal A}$ is defined as
\[
  \Lambda_{P\!I}=\{\lambda\in{\cal A}\,|\,\exists K\in{\cal A}^{(m+n)\times m}\  \lambda T  =K \Delta_I T\}.
\]
We call the ideal $\Lambda_{P\!I}$ the \emph{generalized elementary factor} of
the plant~$P$ with respect to~$I$.  Further, the set of all $\Lambda_{P\!I}$'s
is denoted by~${\cal L}_{P}$, that is, ${\cal L}_{P}=\{\Lambda_{P\!I}\,|\, I\in{\cal I}\}$.
\end{definition}

It is known that the generalized elementary factor of a~plant~$P$ is
independent of the choice of fractions $ND^{-1}=P$ (Lemma\,3.3 of
\cite{bib:mori98bs2}).

The following is the criteria of the feedback stabilizability.
\begin{theorem}\label{Th:3.5}
\textup{\rm (Theorem\,3.2 of \cite{bib:mori98bs2})~} Consider a~causal
plant~$P$.  Then the following statements are equivalent:
\begin{itemize}
\item[\textup{(i)}]
The plant~$P$ is stabilizable.
\item[\textup{(ii)}]
${\cal A}$-modules~${\cal T}_{P,{\cal A}}$ and~${\cal W}_{P,{\cal A}}$ are projective.
\item[\textup{(iii)}]
The set of all generalized elementary factors of~$P$ generates~${\cal A}$;
that is, ${\cal L}_P$ satisfies:
\begin{equation}\label{E:Th:3.5:1}
  \sum_{\Lambda_{P\!I}\in{\cal L}_P} \Lambda_{P\!I}={\cal A}.
\end{equation}
\end{itemize}
\end{theorem}

In the theorem above, each of~(ii) and~(iii) is the necessary and
sufficient conditions of the feedback stabilizability.  Provided that
we can check (\ref{E:Th:3.5:1}) and that we can construct the
right-coprime factorizations over ${\cal A}_{\lambda_I}$ of the given causal
plant, we have given a~method to construct a~causal stabilizing
controller of a~causal stabilizable plant, which has been originally
given in the proof of ``(iii)$\rightarrow$(i)'' of Theorem\,3.2
in~\cite{bib:mori98bs2}.  
We review here the method since it will need later in order to present
a~parameterization of the stabilizing controllers.  

Suppose that (iii) of Theorem\,\ref{Th:3.5} holds.  From
(\ref{E:Th:3.5:1}), there exist~$\lambda_I$'s such that $\sum \lambda_ I=1$,
where~$\lambda_I$ is an element of generalized elementary factor $\Lambda_
{P\!I}$ of the plant~$P$ with respect to~$I\in{\cal I}$; that is, $\lambda_I\in \Lambda_
{P\!I}$.  {Now} let these~$\lambda_ I$'s be fixed.  Further, let ${\cal I}^\sharp$ be
the set of~$I$'s of these~$\lambda_ I$'s, so that
\begin{equation}\label{E:22.Oct.96.205715}
  \sum_{I\in{\cal I}^\sharp} \lambda_I=1.
\end{equation}
We consider without loss of generality that for any $I\in{\cal I}^\sharp$,~$\lambda_I$
is not a~nilpotent element of~${\cal A}$, since $1+x$ is a~unit of~${\cal A}$ for
any nilpotent~$x$ (cf.\,\cite[p.10]{bib:atiyah69a}).  For each $I\in{\cal I}^
\sharp$, ${\cal T}_ {P,{\cal A}_{\lambda_I}}$ is free ${\cal A}_{\lambda_I}$-module of rank~$m$ by
Lemma\,4.10 of \cite{bib:mori98bs2}.  Hence by Lemma\,4.7 of
\cite{bib:mori98bs2}, there exist matrices $\widetilde{X}_I$,
$\widetilde{Y}_I$, $N_I$, $D_I$ over~${\cal A}_ {\lambda_I}$ such that the
following equality holds over~${\cal A}_{\lambda_I}$:
\begin{equation}\label{E:09.Feb.97.151238}
   \widetilde{Y}_IN_I
  +\widetilde{X}_ID_I
  =E_m,
\end{equation}
where $P=N_ID_I^{-1}$ over ${\cal F}({\cal A}_{\lambda_I})$. 

For any positive integer~$\omega$, there are coefficients~$a_I$'s in~${\cal A}$
with $ \sum_{I\in{\cal I}^\sharp} a_I \lambda_I^{\omega}=1$.  We let~$\omega$ be a~sufficiently
large integer.  Hence the matrices $a_I \lambda_I^{\omega}D_I\widetilde{X}_I$ and
$a_I \lambda_I^{\omega} D_I \widetilde{Y}_I$ are over~${\cal A}$ for all $I\in{\cal I}^\sharp$.
Then the causal stabilizing controller has the form
\begin{equation}\label{E:08.Feb.97.215931}
  C=(\sum_{I\in{\cal I}^\sharp} a_I \lambda_I^{\omega} D_I \widetilde{X}_I)^{-1}
    (\sum_{I\in{\cal I}^\sharp} a_I \lambda_I^{\omega} D_I \widetilde{Y}_I). {\ \ } {\ \ }
\end{equation}
In the case where $\sum_{I\in{\cal I}^\sharp} a_I \lambda_I^{\omega} D_I\widetilde{X}_I$ is
${\cal Z}$-singular, we can select other $\widetilde{Y}_I$'s and
$\widetilde{X}_I$'s in order that $\sum_{I\in{\cal I}^\sharp} a_I \lambda_I^{\omega} D_I
\widetilde{X}_I$ is ${\cal Z}$-nonsingular.  For detail, see the proof of
Theorem\,3.2 of~\cite{bib:mori98bs2}.

By using~$C$ of (\ref{E:08.Feb.97.215931}), the matrix $H(P,C)$ in
(\ref{E:H(P,C)}) is calculated as follows: 
\begin{equation}\label{E:31.Jan.99.231643}
   \bmatrix{
    E_n-\sum_{I\in{\cal I}^\sharp} a_I \lambda_I^{\omega} N_I\widetilde{Y}_I
         &  -\sum_{I\in{\cal I}^\sharp} a_I \lambda_I^{\omega} N_I\widetilde{X}_I \cr
    \sum_{I\in{\cal I}^\sharp} a_I \lambda_I^{\omega} D_I \widetilde{Y}_I
         & \sum_{I\in{\cal I}^\sharp} a_I \lambda_I^{\omega} D_I \widetilde{X}_I
  }.
\end{equation}
Since, $\omega$ is a~sufficient large integer, the matrix above is over~${\cal A}$, which implies that the plant~$P$ is stabilized by the
constructed~$C$.

Before finishing this section, we present here a~parallel result of
Lemma\,4.10 of \cite{bib:mori98bs2} for the transposed matrix~$P^t$,
which will be used later, without its proof.  To present it, we
introduce parallel symbols of~${\cal I}$, $I$, and $i_1,\ldots, i_m$.  Let ${\cal J}
$
be the set of all sets of~$n$ distinct integers between~$1$ and $m+n$.
We will use~$J$ as an element of~${\cal J}$ and $j_1,\ldots, j_n$ as elements
of~$J$ with ascending order.
\begin{lemma}\label{L:19.Feb.99.171243}
Consider the transposed matrix~$P^t$.  Let $J\in{\cal J}$.  Let $\Lambda_
{P^t\!J}$ be the generalized elementary factor of the transposed
plant~$P^t$ with respect to~$J$ and further $\sqrt{\Lambda_{P^t\!J}}$
denote the radical of $\Lambda_{P^t\!J}$ (as an ideal).  Suppose that~$\lambda_
J$ is a~non-nilpotent element of $\sqrt{\Lambda_{P^t\!J}}$.  Then, ${\cal A}_{\lambda_
J} $-module ${\cal W}_{P,{\cal A}_{\lambda_J}} $ is free of rank~$n$.
\end{lemma}

\section{Relationship between Generalized Elementary Factors
of the Plant and its Transposed Plant}
\label{S:Stab}
To parameterize (${\cal A}$-)stabilizing controllers, we need to have the
capability to obtain all right-/left-coprime factorizations over~${\cal A}_
{\lambda_I}$, where~$\lambda_I$ is a~nonzero element of the generalized
elementary factor of the plant~$P$ with respect to~$I\in{\cal I}$.  By both
of Lemmas\,4.7 and 4.10 of \cite{bib:mori98bs2}, we already know
that there exist the right-coprime factorization over ${\cal A}_{\lambda_I}$ of
the plant~$P$.  However this is not led to the existence of the
left-coprime factorization over ${\cal A}_{\lambda_I}$ of~$P$.  In this section,
we will show that the radical of any generalized elementary factor of
the plant with respect to~$I\in{\cal I}$, denoted by $\sqrt {\Lambda_ {P\!I}}$,
coincides with that of a~generalized elementary factor of its
transposed transfer matrix with respect to an element~$J\in{\cal J}$, denoted
by $\sqrt {\Lambda_ {P^t\!J}}$; that is, $\sqrt{\Lambda_ {P\!I}}=\sqrt{\Lambda_
{P^t\!J}}$ for appropriate~$I$ and~$J$.  This fact will be led to the
existence of the left-coprime factorizations over ${\cal A}_{\lambda_I}$ of~$P$.

Define first a~bijective mapping~$\tau$ from~${\cal I}$ to~${\cal J}$.  For
convenience we denote by~$I_N$ and~$I_d$ the subsets of~$I$ such that
\[
    I_N =\{ i\,|\, i\leq n, i\in I\},{\ \ }
    I_d =\{ i\,|\, i>n,  i\in I\}.
\]
Using~$I_N$ and~$I_d$, we define~$J_N$ and~$J_d$ as
\[
    J_N = [1,m]\backslash\{ i-n \,|\, i\in I_d\},{\ \ }
    J_d = \{ i+m\,|\, i\in[1,n]\backslash I_N\}.
\]
Then define the bijective mapping $\tau:{\cal I} \rightarrow{\cal J}$ as
\[
  \tau: I_N\cup I_d \mapsto J_N\cup J_d.
\]

Using the mapping~$\tau$, we obtain a~simple result
as follows.
\begin{proposition}\label{Prop:3.1}
$\sqrt{\Lambda_{P\!I}}$ is equal to $\sqrt{\Lambda_{P^t\!\tau(I)}}$.
\end{proposition}
The proof is relatively straightforward and so omitted due to space
limitation.

We now summarize the existence of both of right-/left-coprime
factorizations over~${\cal A}_{\lambda_I}$ of the plant as following.
\begin{proposition}\label{Prop:3:2}
There always exist both of right-/left-coprime factorizations over ${\cal A}_
{\lambda_I}$ of the plant~$P$, where~$\lambda_I$ is a~nonzero element of the
generalized elementary factor of the plant with respect to an~$I$
in~${\cal I}$.
\end{proposition}
\begin{proof}
We already know the existence of the right-coprime factorization over
${\cal A}_{\lambda_I}$ of~$P$.  On the other hand, the existence of the
left-coprime factorization over ${\cal A}_{\lambda_I}$ of~$P$ is derived from
Proposition\,\ref{Prop:3.1}, Lemma\,\ref{L:19.Feb.99.171243} of this
paper and Lemma\,4.7 of \cite{bib:mori98bs2}.  \qquad
\end{proof}

\section{Parameterization of Stabilizing Controllers}
\label{S:01.Feb.99.173244}
In this section, we give a~parameterization of the stabilizing
controllers by one matrix.  

In order to parameterize the stabilizing controllers, we will first
parameterize the transfer matrices $H(P,C)$'s such that~$C$ is
a~stabilizing controller of~$P$ and then obtain the parameterization
of the stabilizing controllers.  We will consider first the case where
there exist both of right-/left-coprime factorizations and then the
general case.

Throughout this section, we assume that~${\cal R}$ denotes~${\cal A}$ or~${\cal A}_f$
and the given plant~$P$ is stabilizable.  For convenience we introduce
two notations used in this section.  Let ${\cal H}(P;{\cal R})$ denote the set of
$H(P,C)$'s such that~$C$ is an ${\cal R}$-stabilizing controller of~$P$ and
${\cal S}(P;{\cal R})$ the set of all ${\cal R}$-stabilizing controllers.  Then the set
${\cal H}(P;{\cal R})$ is expressed as $\{ H(P,C)\,|\, C\in{\cal S}(P;{\cal R})\}$.  Conversely,
once we obtain ${\cal H}(P;{\cal R})$, it is also easy to obtain ${\cal S}(P;{\cal R})$ by the
following lemma.
\begin{lemma}\label{Lemma:4.1}
In this lemma, let $H_{11}\in({\cal R})_n$, $H_{12}\in{\cal R}^{n\times m}$, $H_{21}\in{\cal R}
^{m\times n}$, $H_{22}\in({\cal R})_m$ denote submatrices of~$H$ in $({\cal R})_{m+n}$
such that
\[
   H=\bmatrix{ H_{11} & H_{12} \cr
               H_{21} & H_{22}}.
\]
Then the set of all ${\cal R}$-stabilizing controllers ${\cal S}(P;{\cal R})$ is given
as follows:
\begin{eqnarray}
  {\cal S}(P;{\cal R})
           \!\!&\!\!=\!\!&\!\!
         \lefteqn{\{-\bmatrix{ O_{m\times n}&E_m }H^{-1}
                                   \bmatrix{E_n\cr O_{m\times n}}\,|\,}
                                             \label{E:07}\\
           &&~~~~ H\in{\cal H}(P;{\cal R}), \mbox{~$H$ is nonsingular over~${\cal R}$}\}
                     \nonumber\\
          &&
            \makebox[0mm]{\hspace*{11.4em}$=\{ H_{22}^{-1}H_{21}\,|\, H\in{\cal H}(P;{\cal R}),
             \mbox{~$H$ is nonsingular over~${\cal R}$}\}$}
                     \label{E:08}\\
          &&\makebox[0mm]{\hspace*{11.5em}$=
             \{ H_{21}H_{11}^{-1}\,|\,
                       H\in{\cal H}(P;{\cal R}),\mbox{~$H$ is nonsingular over~${\cal R}$}\}$.}
                     \label{E:09}
\end{eqnarray}
\end{lemma}
\begin{proof}
Since for every $H(P,C)\in{\cal H}(P;{\cal R})$, $\det(E_n+PC)$ is a~nonzerodivisor
of~${\cal R}$, every $H\in{\cal H}(P;{\cal R})$ is nonsingular over~${\cal R}$.  It is known
that the following equality holds:
\[
  H(P,C)=\bmatrix{ E_n & P\cr
                  -C   & E_m}^{-1} .
\]
Directly from this, we have (\ref{E:07}).  The remaining relations
(\ref{E:08}) and (\ref{E:09}) are obtained directly from
(\ref{E:H(P,C)}).  \qquad
\end{proof}

\subsection{The case of existence of right-/left-coprime factorizations}
\label{SS:19.Apr.99.115545}
Throughout this subsection, we assume that the plant~$P$ has both of
right-/left-coprime factorizations over~${\cal R}$.  It should be noted that
from Proposition\,\ref{Prop:3:2}, there exist both of
right-/left-coprime factorizations over ${\cal A}_{\lambda_I}$ of~$P$, where~$\lambda_
I$ is a~nonzero element of the generalized elementary factor of~$P$
with respect to~$I$ in~${\cal I}$.  Hence~${\cal R}$ can be every such ${\cal A}_ {\lambda_
I}$.

Let $(N,D)$ and $(\widetilde{D},\widetilde{N})$ be right-/left-coprime
factorizations over~${\cal R}$ of~$P$, and $\widetilde{Y}_0$,
$\widetilde{X}_0$, $Y_0$, and~$X_0$ be matrices over~${\cal R}$ such that
  $\widetilde{Y}_0N+\widetilde{X}_0D=E_m$ and
  $\widetilde{N}Y_0+\widetilde{D}X_0=E_n$.
We assume here without loss of generality that
\begin{equation}\label{E:22.Apr.99.122126}
  \bmatrix{  \widetilde{X}_0 &  \widetilde{Y}_0 \cr
             \widetilde{N}   & -\widetilde{D}}
  \bmatrix{ D &  Y_0 \cr
            N & -X_0}
  =E_{m+n}.
\end{equation}

The following is a~parameterization of the stabilizing controllers
presented as a~Youla-Ku\v{c}era parameterization.
\begin{theorem}\label{Th:4.2}
\textup{\rm (cf.\,Theorems\,5.2.1 and 8.3.12 of
\cite{bib:vidyasagar85a})~} All matrices~$X$, $Y$, $\widetilde{X}$,
$\widetilde{Y}$ over~${\cal R}$ satisfying
\[
  \widetilde{Y}N+\widetilde{X}D=E_m,{\ \ }
  \widetilde{N}Y+\widetilde{D}X=E_n
\]
are expressed as $X=X_0-NS$, $Y=Y_0+DS$,
$\widetilde{X}=\widetilde{X}_0-R\widetilde{N}$ and
$\widetilde{Y}=\widetilde{Y}_0+R\widetilde{D}$ for~$R$ and~$S$ in ${\cal R}
^{m\times n}$.

Further the set of all ${\cal R}$-stabilizing controllers, denoted by ${\cal S}(P;
{\cal R})$, is given as
\begin{eqnarray*}
  {\cal S}(P;{\cal R})
     \!\!&\!\!=\!\!&\!\!
        \{(\widetilde{X}_0-R\widetilde{N})^{-1}
          (\widetilde{Y}_0+R\widetilde{D})
           \,|\, R\in{\cal R}^{m\times n},\\
     & &\hspace*{5.0em} 
                  \mbox{~$\widetilde{X}_0-R\widetilde{N}$ 
                     is nonsingular over ${\cal R}$}\}\\
     \!\!&\!\!=\!\!&\!\!
        \{(Y_0+DS)(X_0-NS)^{-1}\,|\, S\in{\cal R}^{m\times n},\\
     & &\hspace*{5.0em}
                 \mbox{~$X_0-NS$ is nonsingular over~${\cal R}$}\}.
\end{eqnarray*}
\end{theorem}
This proof is still similar with the previous ones such as
in~\cite{bib:vidyasagar85a}.

The following lemma without the proof gives the set ${\cal H}(P;{\cal R})$
according to Theorem\,\ref{Th:4.2}.
\begin{lemma}\label{Lemma:4.3}
\textup{\rm (cf.\,Corollary\,5.2.7 of \cite{bib:vidyasagar85a})~}
Suppose $P\in{\cal F}({\cal R})^{n\times m}$.  Then ${\cal H}(P;{\cal R})$ is given as
\begin{eqnarray}
&& \makebox[0cm]{$
  {\cal H}(P;{\cal R})=
  \Biggl\{
     \bmatrix{  E_n-N(\widetilde{Y}_0+R\widetilde{D}) 
              &    -N(\widetilde{X}_0-R\widetilde{N}) \cr
                    D(\widetilde{Y}_0+R\widetilde{D}) 
              &     D(\widetilde{X}_0-R\widetilde{N})}
   \Biggm|$\hspace*{2em}}\nonumber\\
&& \makebox[0cm]{$
    R\in{\cal R}^{m\times n}, \mbox{~$\widetilde{X}_0-R\widetilde{N}$ 
                                           is nonsingular over~${\cal R}$}
  \Biggr\};$\hspace*{2em}}
  \label{E:Lemma:4.3:1}
\end{eqnarray}
alternatively
\begin{eqnarray}
&& \makebox[0cm]{$
  {\cal H}(P;{\cal R})=
  \Biggl\{
         \bmatrix{      (X_0-NS)\widetilde{D}
                  &    -(X_0-NS)\widetilde{N}\cr
                        (Y_0+DS)\widetilde{D}
                  & E_m-(Y_0+DS)\widetilde{N}}
   \Biggm| 
$\hspace*{2em}}\nonumber\\
&& \makebox[0cm]{$
    S\in{\cal R}^{m\times n}, \mbox{~$X_0-NS$ is nonsingular over~${\cal R}$}
  \Biggr\}.$\hspace*{2em}}
\label{E:Lemma:4.3:2}
\end{eqnarray}
\end{lemma}

We now start to construct a~new parameterization.  Suppose that~$C_0$
is an ${\cal R}$-stabilizing controller.  We assume without loss of
generality that~$C_0$ is expressed as
$\widetilde{X}_0^{-1}\widetilde{Y}_0$ by Theorem\,\ref{Th:4.2}.
Let~$H_0$ denote $H(P,C_0)\in({\cal R})_ {m+n}$ for short and
$\widehat{H}(R)$ the transfer matrix $H(P,C)$ with
  $C=(\widetilde{X}_0-R\widetilde{N})^{-1}(\widetilde{Y}_0+R\widetilde{D})$.
To include the case where $(\widetilde{X}_0-R\widetilde{N})$ is
singular over~${\cal R}$, we define it as in (\ref{E:Lemma:4.3:1}) as
follows:
\[
  \widehat{H}(R)
  =
         \bmatrix{  E_n-N(\widetilde{Y}_0+R\widetilde{D}) 
                  &    -N(\widetilde{X}_0-R\widetilde{N}) \cr
                        D(\widetilde{Y}_0+R\widetilde{D}) 
                  &     D(\widetilde{X}_0-R\widetilde{N})},
\]
where $R\in{\cal R}^{m\times n}$.  Then the set ${\cal H}(P;{\cal R})$ can be easily expressed
by $\widehat{H}(R)$.
\begin{lemma}\label{Lemma:4.4}
The set ${\cal H}(P;{\cal R})$ is expressed using $\widehat{H}(R)$ as follows:
\begin{eqnarray*}
\lefteqn{
  {\cal H}(P;{\cal R})\!=\!\{\widehat{H}(R)\,|\, R\in{\cal R}^{m\times n},}\\
&&\hspace*{6em} \mbox{~$\widehat{H}(R)$ is nonsingular over~${\cal R}$}\}.
\end{eqnarray*}
\end{lemma}
\begin{proof}
Observe that the following matrix equation holds:
\begin{eqnarray*}
&&  \widehat{H}(R)
  =
   \bmatrix{ E_n & O \cr
             O   & D }
   \bmatrix{  E_n-N(\widetilde{Y}_0+R\widetilde{D}) 
            &    -N \cr
                  \widetilde{Y}_0+R\widetilde{D}
            &     E_m}\times
    \\
&&  \hspace*{12em}
   \bmatrix{ E_n & O \cr
             O   & \widetilde{X}_0-R\widetilde{N}}.
\end{eqnarray*}
The determinants of the matrices in the right hand side of the
equation above are $\det(D)$, $1$,
$\det(\widetilde{X}_0-R\widetilde{N})$, respectively.  Hence the left
hand side is nonsingular over~${\cal R}$ if and only if
$\widetilde{X}_0-R\widetilde{N}$ is nonsingular over~${\cal R}$.  By
(\ref{E:Lemma:4.3:1}) in Lemma\,\ref{Lemma:4.3}, this completes the
proof.  \qquad
\end{proof}

We may define $\widehat{H}'(S)$ denoting $H(P,C)$ with
$C=(Y_0+DS)(X_0-NS)^{-1}$, that is,
\[
  \widehat{H}'(S)
  =
         \bmatrix{      (X_0-NS)\widetilde{D}
                  &    -(X_0-NS)\widetilde{N}\cr
                        (Y_0+DS)\widetilde{D}
                  & E_m-(Y_0+DS)\widetilde{N}},
\]
where $S\in{\cal R}^{m\times n}$ from (\ref{E:Lemma:4.3:2}).  However it can be
easily check $\widehat{H}'(S)=\widehat{H}(S)$ from
(\ref{E:22.Apr.99.122126}).  Hence in the following we are concerned
only with $\widehat{H}(R)$.

We now introduce a~new matrix $\Omega(Q)\in({\cal R})_{m+n}$, which plays a~key
role of new parameterization, as follows:
\begin{equation}\label{E:19.Apr.99.182334}
  \Omega(Q)=
    (
      H_0
      -
      \bmatrix{ E_n \!\!&\!\! O \cr
                O   \!\!&\!\! O}
    )Q
    (
      H_0
      -
      \bmatrix{ O \!\!&\!\! O \cr
                O \!\!&\!\! E_m }
    )
   +
     H_0.
\end{equation}
where $Q\in({\cal R})_{m+n}$.  For convenience define further the partition
of $\Omega(Q)$ as follows:
\[
  \Omega(Q)=\bmatrix{ \Omega_{11}(Q) & \Omega_{12}(Q)\cr
                  \Omega_{21}(Q) & \Omega_{22}(Q)}
\]
with $\Omega_{11}(Q)\in({\cal R})_n$, $\Omega_{12}(Q)\in{\cal R}^{n\times m}$, $\Omega_{21}(Q)\in{\cal R}^{m
\times n}$, $\Omega_{22}(Q)\in({\cal R})_m$.

For the images of the newly introduced matrices, we have the following
relations.
\begin{theorem}\label{Th:4.5}
The images of $\Omega(\cdot)$ and $\widehat{H}(\cdot)$ are
identical; that~is
\begin{eqnarray}
\lefteqn{    \{\Omega(Q)\in({\cal R})_{m+n}\,|\, Q\in({\cal R})_{m+n}\}}\nonumber\\
&&   =       \{\widehat{H}(R)\in({\cal R})_{m+n}\,|\, R\in{\cal R}^{m\times n}\}.
                            \label{E:22.Apr.99.161840:1}
\end{eqnarray}
Further we have
\begin{eqnarray}
\lefteqn{
   {\cal H}(P;{\cal R})=\{\Omega(Q)\in({\cal R})_{m+n}\,|\, Q\in({\cal R})_{m+n},}\nonumber\\
&&\hspace*{6em} \mbox{~$\Omega(Q)$ is nonsingular over~${\cal R}$}\}.
\label{E:22.Apr.99.161840:2}
\end{eqnarray}
\end{theorem}
Here (\ref{E:22.Apr.99.161840:2}) gives a~new parameterization of the
stabilizing controllers.
\begin{proof}
We prove only (\ref{E:22.Apr.99.161840:1}) by showing that (i) for any
matrix~$Q$, there exists a~matrix~$R$ such that $ \widehat{H}(R)=\Omega
(Q)$, and that (ii) for any matrix~$R$, there exists a~matrix~$Q$ such
that $\widehat{H}(R)=\Omega(Q)$.  Once (\ref{E:22.Apr.99.161840:1}) is
obtained, (\ref{E:22.Apr.99.161840:2}) is obvious by
Lemma\,\ref{Lemma:4.4}.

We first prove~(i).  We rewrite (\ref{E:19.Apr.99.182334}) as follows:
\begin{eqnarray}
&&\hspace*{3em}
   \makebox[0cm]{$  \Omega(Q)=
    \bmatrix{ O & -N \cr
              O &  D }
    \bmatrix{ X_0 & N \cr
             -Y_0 & D }^{-1}
    Q
    \bmatrix{ -\widetilde{Y}_0 & \widetilde{X}_0 \cr
               \widetilde{D}   & \widetilde{N}}^{-1}\times$}\nonumber\\
&&    \bmatrix{  O & O \cr
              \widetilde{D} & -\widetilde{N}}
    +
    H_0.\label{E:16.Apr.99.133905}
\end{eqnarray}
Since both of the inverse matrices in (\ref{E:16.Apr.99.133905}) are
unimodular (cf.~Corollary\,4.1.67 of~\cite{bib:vidyasagar85a}), we let
\[
  Q'
  =
  \bmatrix{ X_0 & N \cr
           -Y_0 & D }^{-1}
  Q
    \bmatrix{ -\widetilde{Y}_0 & \widetilde{X}_0 \cr
               \widetilde{D}   & \widetilde{N}}^{-1}.
\]
Then (\ref{E:16.Apr.99.133905}) can be rewritten further as follows:
\begin{equation}\label{E:15.Apr.99.230806:x:y}
  \Omega(Q)=
  \bmatrix{ O & -N \cr
            O &  D }
  Q'
  \bmatrix{  O & O \cr
            \widetilde{D} & -\widetilde{N}}
  +
   H_0.
\end{equation}
Partition~$Q'$ as
\[
  \bmatrix{ Q_{11}' & Q_{12}' \cr
            Q_{21}' & Q_{22}'}
  =
  Q',
\]
where $Q_{11}'\in{\cal R}^{n\times m}$, 
      $Q_{12}'\in({\cal R})_n$, 
      $Q_{21}'\in({\cal R})_m$, 
      $Q_{22}'\in{\cal R}^{m\times n}$. 
Then (\ref{E:15.Apr.99.230806:x:y}) can be rewritten again as follows:
\[
  \Omega(Q)=
      \bmatrix{  E_n-N(\widetilde{Y}_0+Q_{22}'\widetilde{D}) 
               &    -N(\widetilde{X}_0-Q_{22}'\widetilde{N}) \cr
                     D(\widetilde{Y}_0+Q_{22}'\widetilde{D}) 
               &     D(\widetilde{X}_0-Q_{22}'\widetilde{N})},
\]
which is equal to $\widehat{H}(Q_{22}')$.  Therefore the matrix
$Q_{22}'$ is the matrix~$R$ satisfying $\widehat{H}(R)=\Omega(Q)$.

Next we prove (ii).  From the proof of (i), letting~$Q$ as
\[
   Q=
     \bmatrix{ X_0 & N \cr
              -Y_0 & D }
     \bmatrix{ \times & \times \cr
               \times & R}
     \bmatrix{ -\widetilde{Y}_0 & \widetilde{X}_0 \cr
                \widetilde{D}   & \widetilde{N}},
\]
where $\times$'s denote arbitrary matrices, we obtain directly
$\widehat{H}(R)=\Omega(Q)$.  \qquad
\end{proof}

\subsection{The General Case}
\label{SS:25.Apr.99.212011}
In this subsection, we parameterize all stabilizing controllers over~$
{\cal A}$ even in the case where there does not exist right-/left-coprime
factorizations over~${\cal A}$.

Let~$\lambda_I$ denote an arbitrary but fixed element of the generalized
elementary factor of~$P$ with respect to~$I\in{\cal I}^\sharp$ satisfying
(\ref{E:22.Oct.96.205715}).  As mentioned
in~\S\,\ref{SS:19.Apr.99.115545}, there exist both of
right-/left-coprime factorizations over ${\cal A}_{\lambda_I}$ of the plant~$P$.
We let $(N_I,D_I)$ and $(\widetilde{D}_I,\widetilde{N}_I)$ denote
right-/left-coprime factorizations over ${\cal A}_{\lambda_I}$ of~$P$.  Suppose
again that~$C_0$ is a~stabilizing controller.  By
Theorem\,\ref{Th:4.2}, there exist right-/left-coprime factorizations
$(Y_{0I},X_{0I})$ and $(\widetilde{X}_{0I},\widetilde{Y}_{0I})$ over $
{\cal A}_{\lambda_I}$ of~$C_0$ such that
  $\widetilde{Y}_{0I}N_I+\widetilde{X}_{0I}D_I=E_m$ and
  $\widetilde{N}_IY_{0I}+\widetilde{D}_IX_{0I}=E_n$.
In order to distinguish~$\widehat{H}$ for each $I\in{\cal I}^\sharp$, we
introduce~$\widehat{H}_I$ analogously to~\S\,\ref{SS:19.Apr.99.115545}
as follows:
\[
  \widehat{H}_I(R_I)\!=\!
         \bmatrix{  E_n-N_I(\widetilde{Y}_{0I}+R_I\widetilde{D}_I) 
                  &    -N_I(\widetilde{X}_{0I}-R_I\widetilde{N}_I) \cr
                        D_I(\widetilde{Y}_{0I}+R_I\widetilde{D}_I) 
                  &     D_I(\widetilde{X}_{0I}-R_I\widetilde{N}_I)},
\]
where~$R_I\in{\cal A}_{\lambda_I}^{m\times n}$.  On the other hand, the matrices $\Omega
(Q)$ and $\Omega_ {ij}(Q)$'s with $i,j=1,2$ are still used.

In the following we show that (\ref{E:22.Apr.99.161840:2}) in
Theorem\,\ref{Th:4.5} still holds even in the case where there do not
exist right-/left-coprime factorizations over~${\cal A}$.
\begin{theorem}\label{Th:4.6}
The following equality holds:
\begin{eqnarray*}
&& {\cal H}(P;{\cal A})=\{\Omega(Q)\in({\cal A})_{m+n}\,|\, Q\in({\cal A})_{m+n},\\
&& \hspace*{13em}\mbox{~$\Omega(Q)$ is nonsingular}\}.
\end{eqnarray*}
\end{theorem}
\begin{proof}
In order to prove this theorem, it is sufficient to show the
following:
\begin{itemize}
\item[(i)] For any matrix~$Q$, if $\Omega(Q)$ is nonsingular, there
exists a~stabilizing controller~$C$ such that $H(P,C)=\Omega(Q)$.
\item[(ii)] Conversely, for any stabilizing controller~$C$, there exists
a~matrix~$Q$ such that $H(P,C)=\Omega(Q)$.
\end{itemize}

We first prove (i).  Suppose that $\Omega(Q)$ is nonsingular.  Assume
without loss of generality that $\lambda_{I_0}$ is a~nonzerodivisor with
$I_0\in{\cal I}^\sharp$.  Then by Theorem\,\ref{Th:4.5}, there exists
a~matrix~$R_{I_0}$ over ${\cal A}_{\lambda_{I_0}}$ such that
$\widehat{H}_{I_0}(R_{I_0})=\Omega(Q)$.  By Lemma\,\ref{Lemma:4.4}, there
exists an ${\cal A}_ {\lambda_ {I_0}}$-stabilizing controller~$C$.  Observe now
that ${\cal F}({\cal A}_{\lambda_ {I_0}})={\cal F}$, so that~$C$ is over~${\cal F}$.  Since
$H(P,C)=\Omega(Q)$, $H(P,C)$ is over~${\cal A}$, which implies that~$C$ is
a~stabilizing controller of~$P$.

Next we prove (ii).  Suppose that~$P$ is stabilizable.  As
in~\S\,\ref{S:PR}, let $a_I$'s be in~${\cal A}$ for $I\in{\cal I}^\sharp$ such that $
\sum_ {I\in{\cal I}^\sharp} a_I \lambda_I^{\omega}=1$ with a~sufficiently large integer~$\omega
$.  Let~$C$ be an arbitrary but fixed stabilizing controller of~$P$.
Since~$C$ is also an ${\cal A}_{\lambda_I}$-stabilizing controller, by
Lemma\,\ref{Lemma:4.4} there exists~$R_I$ over ${\cal A}_{\lambda_I}$ such that
  $H(P,C)=\widehat{H}_I(R_I)$ 
for each $I\in{\cal I}^\sharp$.  By Theorem\,\ref{Th:4.5} there
exists a~matrix~$Q_I$ over ${\cal A}_{\lambda_I}$ such that $ \widehat{H}_I(R_I)=
\Omega(Q_I)$ for each $I\in{\cal I}^\sharp$, so that
  $H(P,C)=\Omega(Q_I)$ 
over ${\cal A}_{\lambda_I}$.  Hence we have
\[
  H(P,C)=\sum_{I\in{\cal I}^\sharp} a_I \lambda_I^{\omega}\Omega(Q_I)
\]
over~${\cal A}$.  From (\ref{E:19.Apr.99.182334}), the equation above can be
rewritten as follows:
\begin{eqnarray*}
  && H(P,C)\\
        \!&\!=\!&\!
         \sum_{I\in{\cal I}^\sharp} a_I \lambda_I^{\omega}
                              (
                                H_0
                                -
                                \bmatrix{ E_n      \!\!&\!\! O \cr
                                          O \!\!&\!\! O}
                              )Q_I
                              (
                                H_0
                                -
                                \bmatrix{ O \!\!&\!\! O \cr
                                          O \!\!&\!\! E_m }
                              )
                             +
                             H_0 \\
        \!&\!=\!&\!\Omega(\sum_{I\in{\cal I}^\sharp} a_I \lambda_I^{\omega}Q_I).
\end{eqnarray*}
Letting $Q=\sum_{I\in{\cal I}^\sharp} a_I \lambda_I^{\omega}Q_I$, we have proved (ii).
\qquad
\end{proof}

We now have a~parameterization of the stabilizing controllers over~${\cal A}
$ by virtue of Lemma\,\ref{Lemma:4.1}.
\begin{corollary}\label{Cor:4.7}
The set of all stabilizing controllers ${\cal S}(P;{\cal A})$ is given as follows:
\begin{eqnarray*}
  {\cal S}(P;{\cal A})
           &=&
         \{-\bmatrix{ O_{m\times n}&E_m }\Omega(Q)^{-1}
                                  \bmatrix{E_n\cr O_{m\times n}}\,|\,\\
           &&~~~~~~~~~ Q\in({\cal A})_{m+n}, \mbox{~$\Omega(Q)$ is nonsingular}\}\\
&& \makebox[0pt]{\hspace*{11em}$=\{\Omega_{22}(Q)^{-1}\Omega_{21}(Q)\,|\,
                       Q\in({\cal A})_{m+n},\mbox{~$\Omega(Q)$ is nonsingular}\}$}\\
&& \makebox[0pt]{\hspace*{11.1em}$=\{\Omega_{21}(Q)\Omega_{11}(Q)^{-1}\,|\,
                       Q\in({\cal A})_{m+n},\mbox{~$\Omega(Q)$ is nonsingular}\}$.}
\end{eqnarray*}
\end{corollary}

To finish this section, we review the causality of the stabilizing
controllers of strictly causal plants.
\begin{proposition}\label{Prop:4.8}
\textup{(Proposition\,6.2 of~\cite{bib:mori98bs2})~} For any
stabilizable strictly causal plant, all stabilizing controllers of the
plant must be causal.
\end{proposition}

As a~result, when the plant is strictly causal, all stabilizing
controllers obtained by the parameterization of
Corollary\,\ref{Cor:4.7} are causal.

\section{Examples of Parameterization} 
\label{S:Example}
Let us consider here the multidimensional systems with structural
stability\cite{bib:guiver85a}.  Recall that the construction methods
of a~stabilizing controller have already been presented by
Sule\cite{bib:sule94a} and Lin\cite{bib:lin98b}.  Thus, as an example
of the parameterization, using Corollary\,\ref{Cor:4.7} we can obtain
the parameterization of the stabilizing controllers.  This gives one of
the solutions of an open problem about the parameterization of the
stabilizing controllers of the multidimensional systems given by
Lin\cite{bib:lin98a}.

The reader can further refer to~\cite{bib:mori99bta} for other
examples of parameterizations of the stabilizing controllers.

\section{Concluding Remarks and Further Works}
In this paper we have given a~parameterization of the stabilizing
controllers over commutative rings.  In this paper the minimal number
of parameters required to give the parameterization is not clarified.
Since it is important for the implementation of the parameterization,
it should be clarified as a~further problem.

%
%


\end{document}